\documentclass[11pt]{amsart}
\usepackage{amsfonts}
\usepackage{graphicx,color}
\usepackage{amsthm}
\usepackage{amssymb,amsmath}

\title{\bf On spherical codes with inner products in a prescribed interval}

\date{\today}

\newtheorem{theorem}{Theorem}[section]
\newtheorem{lemma}[theorem]{Lemma}
\newtheorem{corollary}[theorem]{Corollary}

\newtheorem{conjecture}[theorem]{Conjecture}

\theoremstyle{definition}
\newtheorem{definition}[theorem]{Definition}

\newtheorem{remark}[theorem]{Remark}

\author[P. Boyvalenkov]{P. G. Boyvalenkov $^\dagger$}
\address{Institute of Mathematics and Informatics, Bulgarian Academy of Sciences,
8 G Bonchev Str.,
1113  Sofia, Bulgaria \\
and Faculty of Mathematics and Natural Sciences, South-Western University, Blagoevgrad, Bulgaria.
}
\email{peter@math.bas.bg}
\thanks{\noindent $^\dagger$ The research of these authors was supported, in part, by a Bulgarian NSF contract DN02/2-2016.}

\author[P. Dragnev]{P. D. Dragnev $^{\dagger \dagger}$}
\address{Department of Mathematical Sciences,
Purdue University
Fort Wayne, IN 46805, USA }
\email{dragnevp@ipfw.edu}
\thanks{\noindent $^{\dagger \dagger}$ The research of this author was supported, in part, by a Simons Foundation grant no. 282207.}

\author[D. Hardin]{D. P. Hardin$^*$}
\address{Center for Constructive Approximation, Department of Mathematics, \hspace*{.1in}
Vanderbilt University,
Nashville, TN 37240, USA  }
\email{doug.hardin@vanderbilt.edu}

\author[E. Saff]{E. B. Saff$^*$}
\email{edward.b.saff@vanderbilt.edu}
\thanks{\noindent $^*$ The research of these authors was supported, in part,
by the U. S. National Science Foundation under grant DMS-1516400.
}
\author[M. Stoyanova]{M. M. Stoyanova$^\dagger$}
\address{Faculty of Mathematics and Informatics,
Sofia University,
5 James Bourchier Blvd.,
1164 Sofia, Bulgaria}
\email{stoyanova@fmi.uni-sofia.bg}

\begin{document}
\maketitle

\begin{abstract}
We develop a framework for obtaining linear programming bounds for spherical
codes whose inner products belong to a prescribed subinterval $[\ell,s]$ of $[-1,1)$.
An intricate relationship between Levenshtein-type upper bounds on cardinality of codes
with inner products in $[\ell,s]$ and lower bounds on the potential energy (for absolutely monotone interactions)
for codes with inner products in $[\ell,1)$ (when the cardinality of the code is kept fixed)
is revealed and explained. Thereby, we obtain a new extension of Levenshtein bounds for such codes.
The universality of our bounds is exhibited by a unified derivation and their validity for a wide range of codes and potential functions.
\end{abstract}



\section{Introduction}

In the seminal paper of Cohn and Kumar \cite{CK}, many classical maximal spherical codes
with applications to communications, such as the Korkin-Zolotarev kissing number configuration on $\mathbb{S}^7$,
the Leech lattice configuration in 24 dimensions, the $600$-cell, etc.,
were shown to be universally optimal in the sense that they have minimal potential energy for a large class of potential interactions.
The notion of universal optimality was further developed for Hamming spaces in \cite{CZ}.

As important as these particular configurations are, it is of significant interest to study bounds for codes of general cardinality. The theory of universal bounds for codes and designs in polynomial metric spaces was laid out by Levenshtein in \cite{Lev}. The interplay between Levenshtein's framework
and universal lower bounds (ULB) on potential energy of codes was established recently by the authors for Euclidean spaces in \cite{BDHSS-CA} and for
Hamming spaces in \cite{BDHSS-DCC}. In this paper we further that interplay to codes with inner products
in a prescribed subinterval $[\ell,s]$ of $[-1,1)$ and as a result derive an extension of Levenshtein's framework to this setting.

Let $\mathbb{S}^{n-1}\subset \mathbb{R}^n$ denote the $(n-1)$-dimensional unit sphere. A nonempty finite set $C \subset \mathbb{S}^{n-1}$ is called a {\em spherical code}.
For $-1 \leq \ell <s < 1$ denote by
\[ \mathcal{C}(\ell,s):=\{C \subset \mathbb{S}^{n-1}: \ \ell \leq \langle x,y \rangle \leq s , \ x,y \in C, x \neq y\}, \]
the set of spherical codes with prescribed maximum diameter and minimum pairwise distance, where
$\langle x,y \rangle$ denotes the inner product of $x$ and $y$.
We establish upper bounds on the quantity
\[ \mathcal{A}(n;[\ell,s]):=\max\{|C|:C \in \mathcal{C}(\ell,s)\}, \]
which is a classical problem in coding theory.


Given a (potential) function $h(t):[-1,1] \to [0,+\infty]$ and a code $C \subset \mathbb{S}^{n-1}$, we define the {\em potential energy} (also referred to as {\em $h$-energy}) of $C$ as
\[ E(C;h):=\sum_{x, y \in C, x \neq y} h(\langle x,y \rangle). \]
In what follows we shall consider potential functions $h$ that are absolutely monotone, namely $h^{(k)}(t) \geq 0$ for every $k \geq 0$ and $t\in [-1,1)$. For such potentials we establish ULB for the quantity
\[ \mathcal{E}(n,M,\ell;h):=\inf\{E(C;h):C \in \mathcal{C}(\ell,1), |C|=M \}. \]
As in \cite{BDHSS-CA}, the use of linear programming reveals a strong connection between our
ULB on $\mathcal{E}(n,M,\ell;h)$ and our Levenshtein-type upper bounds on $\mathcal{A}(n;[\ell,s])$.

Throughout, $P_k^{(n)}(t)$, $k=0,1,\dots$, will denote the Gegenbauer polynomials \cite{Sze}
normalized with $P_k (1)=1$.
We consider functions $f(t):[-1,1]\to \mathbb{R}$,
\[ f(t)=\sum_{k=0}^\infty f_k P_k^{(n)}(t), \ \ {\rm where} \ \
f(1)=\sum_{k=0}^\infty f_k <\infty.\]
The function $f$ is called {\em positive definite} (\emph{strictly positive definite}) if all coefficients $f_k$ are non-negative (positive). Following Levenshtein's notation we denote the class of all positive definite  (strictly positive definite) functions by the symbol ${\mathcal F}_{\geq}$ (${\mathcal F}_{>}$). When $f$ is a polynomial, the definition of ${\mathcal F}_{>}$ does not include $f_k$ for $k>\deg(f)$ (since $f_k=0$ for such $k$).

The Kabatiansky-Levenshtein \cite{KL} approach (see also \cite{DGS}) is based on the inequality
\begin{equation}\label{MaxCodesLP} \mathcal{A}(n;[\ell,s]) \leq \min_{f\in \mathcal{F}_{n,\ell,s}} f(1)/f_0,\end{equation}
where
\[ \mathcal{F}_{n,\ell,s}:=\{ f \in {\mathcal F}_{\geq} \ | \ f(t) \leq 0, \,  t \in [\ell,s], \, f_0>0\}. \]

Similarly, the Delsarte-Yudin approach (see \cite{Y}) uses the inequality
\begin{equation}\label{EnergyLP} \mathcal{E}(n,M,\ell;h)\ge \max_{g\in \mathcal{G}_{n,\ell;h}} M(Mg_0-g(1)), \end{equation}
where
\[ \mathcal{G}_{n,\ell;h}:=\{ g\in {\mathcal F}_{\geq} \ | \  g(t) \leq h(t), \, t \in [\ell,1),  \, g_0>0\}. \]

The determination of the right-hand sides of the bounds \eqref{MaxCodesLP} and \eqref{EnergyLP} over the respective classes defines two infinite linear programs. To determine his universal bounds on $\mathcal{A}(n,s):=\mathcal{A}(n;[-1,s])$ Levenshtein \cite{Lev} found explicitly the solution of the linear program posed by \eqref{MaxCodesLP} when restricted to $\mathcal{F}_{n,-1,s}\cap \mathcal{P}_m$, where $\mathcal{P}_m$ denotes the class of real polynomials of degree at most $m$.

In \cite{BDHSS-CA} the authors considered the linear program in \eqref{EnergyLP} over $\mathcal{G}_{n,-1;h} \cap \mathcal{P}_m$ and found its solution as the Hermite interpolation polynomial of $h(t)$ at the zeros of the Levenshtein polynomial. This implies the ULB on $\mathcal{E}(n,M;h):=\mathcal{E}(n,M,-1;h)$. The interplay between the two optimal solutions is that the zeros of the Levenshtein polynomials serve also as nodes of an important Radau or Lobato quadrature formulae.

In this paper we further develop the intricate connection
between the maximum cardinality and minimum energy problems, which is described in our main result Theorem \ref{main-th-subintervals}. For this purpose a central role is played by an $\ell$-modification of the so-called `strengthened Krein condition' introduced by Levenshtein (see Section 4).

The outline of the paper is as follows. In Section 2, we introduce certain signed measures and establish their positive definiteness up to an appropriate degree. Properties of their associated orthogonal polynomials are also discussed. In Section 3, Levenshtein-type polynomials $f_{2k}^{(n,\ell,s)}(t)$ are constructed and corresponding quadrature formulas are derived. These formulas are used in Section 4, together with linear programming techniques, to derive the Levenshtein-type bounds on the cardinality of maximal codes $\mathcal{A}(n;[\ell,s])$ and ULB-type
(in the sense of \cite{BDHSS-CA}) energy bounds on $\mathcal{E}(n,M,\ell;h)$. In the last
section
some special examples and numerical evidence of an $\ell$-strengthened Krein property are presented.

%
\section{Positive definite signed measures and associated orthogonal polynomials}

In this section we establish the positive definiteness up to certain degrees
of the signed measures that are used in the proof of our main result, Theorem \ref{main-th-subintervals}.

We shall denote the measure of orthogonality of Gegenbauer polynomials as
\begin{equation}\label{dmu} d\mu(t):=\gamma_n (1-t^2)^{\frac{n-3}{2}}\, dt, \quad t\in [-1,1], \quad
 \gamma_n := \frac{\Gamma(\frac{n}{2})}{\sqrt{\pi}\Gamma(\frac{n-1}{2})}, \end{equation}
where $ \gamma_n$ is a normalizing constant that makes $\mu$ a probability measure.

Levenshtein used the {\em adjacent} (to Gegenbauer) {\em polynomials}
\begin{equation}\label{Adjacent}
P_k^{1,0} (t): =P_k^{(\frac{n-1}{2},\frac{n-3}{2})}(t)/P_k^{(\frac{n-1}{2},\frac{n-3}{2})}(1)=
\eta_k^{1,0} t^k + \cdots ,\quad \eta_k^{1,0}>0,
\end{equation}
where $P_k^{(\alpha,\beta)}(t)$ denotes the classical Jacobi polynomial (the normalization is again chosen so that $P_k^{1,0} (1)=1$).  The polynomials \eqref{Adjacent} are orthogonal with respect to the probability measure
\begin{equation}\label{AdjacentMeas} d\chi(t) := (1-t) d\mu(t).\end{equation}
They also satisfy the following three-term recurrence relation
\[ (t-a_i^{1,0})P_i^{1,0} (t)=b_i^{1,0} P_{i+1}^{1,0} (t)+c_i^{1,0} P_{i-1}^{1,0} (t), \quad i=1,2, \dots , \]
where
\[ P_0^{1,0} (t) = 1, \ P_1^{1,0} (t) = \frac{nt+1}{n+1}, \]
\[ b_i^{1,0} =\frac{\eta_i^{1,0}}{\eta_{i+1}^{1,0}}>0, \
c_i^{1,0} =\frac{r_{i-1}^{1,0} b_{i-1}^{1,0}}{r_i^{1,0}}, \ a_i^{1,0}=1-b_i^{1,0}-c_i^{1,0}, \]
\[r_i^{1,0}:=\left( \int_{-1}^1 \left[ P_i^{1,0} (t)\right]^2 d\chi(t) \right)^{-1} = \left(\frac{n+2i-1}{n-1}\right)^2\binom{n+i-2}{i}.\]

\noindent Let $t_{i,1}^{1,0}< t_{i,2}^{1,0} < \dots < t_{i,i}^{1,0}$ be the zeros of the polynomial $P_i^{1,0} (t)$, which are known to interlace with the zeros of $P_{i-1}^{1,0} (t)$.

We next recall the definition of positive definite signed measures up to degree $m$ (see \cite[Definition 3.4]{CK}).

\begin{definition}
A signed Borel measure $\nu$ on $\mathbb{R}$ for which all polynomials are integrable
is called {\bf positive definite up to degree $\mathbf {m}$} if for all real polynomials $p \not\equiv 0$
of degree at most $m$ we have $\int p(t)^2 d \nu(t) > 0$.
\end{definition}

Given $\ell$ and $s$ such that $\ell < t_{k,1}^{1,0}<t_{k,k}^{1,0}<s$, we define the signed measures on $[-1,1]$ (see \eqref{dmu} and \eqref{AdjacentMeas})
\begin{eqnarray}
\label{SignedMeasures1}
d\nu_\ell(t) &:=& (t-\ell)d\chi(t),\\
\label{SignedMeasures2}
d\nu_s(t) &:=& (s-t)d\chi(t), \\
\label{SignedMeasures3}
d\nu_{\ell,s}(t) &:=& (t-\ell)(s-t)d\chi(t), \\
\label{SignedMeasures4}
d\mu_\ell (t) &:=& (t-\ell )d\mu (t).
\end{eqnarray}
The following lemma establishes the positive definiteness of these signed measures
up to certain degrees, which in turn allows us to define orthogonal polynomials
with respect to these signed measures. This equips us with the essential ingredients
for modifying Levenshtein's framework.

\begin{lemma} \label{lem_pos_def} For given $k>1$, let $s$ and $\ell$  satisfy $\ell < t_{k,1}^{1,0}<t_{k,k}^{1,0}<s$.
Then the measures $d\nu_\ell(t)$, $d\nu_s (t)$, and $d\mu_\ell (t)$ are
positive definite up to degree $k-1$ and the measure $d\nu_{\ell,s} (t)$ is
positive definite up to degree $k-2$.
\end{lemma}

{\it Proof.} We first note that the system of $k+1$ nodes
\[M_{k+1}:=\{ t_{k,1}^{1,0}< t_{k,2}^{1,0}< \dots< t_{k,k}^{1,0}<1:=t_{k,k+1}^{1,0} \}\]
defines a positive Radau quadrature with respect to the measure $\mu$
that is exact for all polynomials of degree at most $2k$ (see e.g. \cite[pp. 102-105]{DR}, \cite[Theorem 2.4]{BCV}), namely the quadrature formula
\begin{equation}\label{w_quadrature}
f_0:= \int_{-1}^1 f(t) d\mu(t)= w_{k+1} f(1)+\sum_{i=1}^k w_i f(t_{k,i}^{1,0} )
\end{equation}
holds for all polynomials $f$ of degree at most $2k$, and the weights $w_i$, $i=1,
\dots, k+~1$ are positive.
%
%
%
%

Let now $q(t)$ be an arbitrary polynomial of degree at most $k-1$. From \eqref{w_quadrature} we have that
\begin{eqnarray*}
\int_{-1}^1 q^2(t) d\nu_\ell (t) &=& \int_{-1}^1 q^2(t)(1-t)(t-\ell) d\mu(t) \\
&=& \sum_{i=1}^k w_i q^2(t_{k,i}^{1,0} )(1-t_{k,i}^{1,0} )(t_{k,i}^{1,0} -\ell) \geq 0,
\end{eqnarray*}
where equality may hold only if $q(t_{k,i}^{1,0} )=0$ for all $i=1,\dots,k$,
which would imply that $q(t) \equiv 0$. Therefore the measure $d\nu_\ell (t)$ is positive definite up to degree $k-1$
as asserted.

Similarly, for the measure $d\nu_s (t)$ and ${\rm deg}\ q \leq k-1$ we have
\begin{eqnarray*}
\int_{-1}^1 q^2(t) d\nu_s (t) &=& \int_{-1}^1 q^2(t)(1-t)(s-t) d\mu(t) \\
&=& \sum_{i=1}^k w_i q^2(t_{k,i}^{1,0} )(1-t_{k,i}^{1,0} )(s-t_{k,i}^{1,0} ) \geq 0,
\end{eqnarray*}
where again equality holds only for $q(t) \equiv 0$.

Next, if $q(t) \not\equiv 0$ is of degree at most $k-2$, then
we utilize \eqref{w_quadrature} again to derive that
\[
 \int_{-1}^1 q^2(t) d\nu_{\ell,s} (t) = \sum_{i=1}^k w_i q^2(t_{k,i}^{1,0} )(1-t_{k,i}^{1,0} )(t_{k,i}^{1,0} -\ell)(s-t_{k,i}^{1,0} ) > 0.
\]
Hence, $d\nu_{\ell,s} (t)$ is positive definite up to degree $k-2$.

To verify the assertion about the measure $d\mu_\ell (t)$ we employ a similar argument
but with a quadrature rule defined on the collection of $k$ nodes
\[\widetilde{M}_k:=\{ t_{k,1}< t_{k,2}< \cdots < t_{k,k} \},\]
where $t_{k,i}$ are the zeros of the regular Gegenbauer polynomials $P_k^{(n)}(t)$.
We note that from \cite[Lemma 5.29, Eq. (72)]{Lev} we have $t_{k,1}^{1,0}<t_{k,1}$.
Using the associated Lagrange basis polynomials $\widetilde{L}_i$, $i=1,2,\dots, k$, we define
the weights $v_i:= \int_{-1}^1 \widetilde{L}_i (t) d\mu(t)$, $i=1,2, \dots, k$.
Then, as in the proof of Gaussian quadrature, one shows that the formula
\[ f_0:= \int_{-1}^1 f(t) d\mu(t)= \sum_{i=1}^k v_i f(t_{k,i} ) \]
is exact for polynomials of degree up to $2k-1$.
Thus, for any polynomial $q(t)$ of degree less than or equal to $k-1$, we have
\begin{eqnarray*}
\int_{-1}^1 q^2(t) d\mu_\ell (t) &=& \int_{-1}^1 q^2(t)(t-\ell) d\mu(t) \\
&=& \sum_{i=1}^k v_i q^2(t_{k,i} )(t_{k,i} -\ell) \geq 0,
\end{eqnarray*}
with equality if and only if $q(t)\equiv 0$. This concludes the proof of the lemma. \hfill $\Box$

\medskip

Applying Gram-Schmidt orthogonalization (see, for example, \cite[Lemma 3.5]{CK}) one derives the
existence and uniqueness (for the so-chosen normalization) of the following classes
of orthogonal polynomials with respect to the signed measures \eqref{SignedMeasures1}-\eqref{SignedMeasures4}.

\begin{corollary} \label{cor_ortho} Let $\ell < t_{1,k}^{1,0}<t_{k,k}^{1,0}<s$.
The following classes of orthogonal polynomials are well-defined:
\[ \{ P_j^{0,\ell} (t) \}_{j=0}^k, \ {\rm orthogonal\ w.r.t.} \ d\mu_\ell (t), \ P_j^{0,\ell} (1)=1; \]
\[ \{ P_j^{1,\ell} (t) \}_{j=0}^k, \ {\rm orthogonal\ w.r.t.} \ d\nu_\ell (t), \ P_j^{1,\ell} (1)=1; \]
\[ \{ P_j^{1,s} (t) \}_{j=0}^k, \ {\rm orthogonal\ w.r.t.} \ d\nu_s (t),  \ P_j^{1,s} (1)=1; \]
\[ \{ P_j^{1,\ell,s} (t) \}_{j=0}^{k-1}, \ {\rm orthogonal\ w.r.t.}\ d\nu_{\ell,s} (t),  \ P_j^{1,\ell,s} (1)=1. \]
The polynomials in each class satisfy a three-term recurrence relation and their zeros interlace.
\end{corollary}

\begin{remark} We note that if $t_{k+1,1}^{1,0}<\ell<t_{k,1}^{1,0}$ is such that
$P_{k+1}^{1,0}(\ell)/P_k^{1,0}(\ell)=1$, then $P_k^{1,\ell}(1)=0$ and the
normalization above fails. However, for our purposes we shall restrict to values
of $\ell$ such that $P_{k+1}^{1,0}(\ell)/P_k^{1,0}(\ell)<1$.
\end{remark}

Utilizing the Christoffel-Darboux formula  (see, for example \cite[Th. 3.2.2]{Sze}, \cite[Eq. (5.65)]{Lev})
we are able to construct these polynomials explicitly.
Let
\begin{equation}
\label{kernelT}
\begin{split}
T_i^{1,0} (x,y) &:= \sum_{j=0}^i r_j^{1,0} P_j^{1,0}(x)  P_j^{1,0}(y) \\
&= r_i ^{1,0}b_i^{1,0} \frac{P_{i+1}^{1,0}(x)  P_i^{1,0}(y)-P_{i+1}^{1,0}(y)  P_i^{1,0}(x)}{x-y}.
\end{split}
\end{equation}
Note that in the limiting case $x=y$ we use appropriate derivatives.

Levenshtein \cite{Lev} uses the Christoffel-Darboux formula to prove the
interlacing properties $t_{j+1,i}<t_{j,i}^{1,0}<t_{j,i}$, $i=1,2,\dots,j$,
of the zeros of $P_j^{1,0}$ and the Gegenbauer polynomials. Similarly, from
the representation
\[ P_j^{1,\ell} (t)=\frac{(1-\ell)\left( P_{j+1}^{1,0} (t)-P_j^{1,0} (t)P_{j+1}^{1,0} (\ell)/P_{j}^{1,0} (\ell) \right)}{(t-\ell)\left(1-P_{j+1}^{1,0} (\ell)/P_{j}^{1,0} (\ell)\right)},\]
which is verified in the next theorem, we derive interlacing properties of the zeros of $P_j^{1,\ell}$ with respect to the zeros of $P_i^{1,0}$.
%
%

\medskip

\begin{theorem}
\label{rootsofP1l}
Let $\ell$ and $k$ be such that $t_{k+1,1}^{1,0}<\ell<t_{k,1}^{1,0}$ and
$P_{k+1}^{1,0}(\ell)/P_{k}^{1,0}(\ell)<1$. Then all zeros $\{ t_{i,j}^{1,\ell}\}_{j=1}^i$ of $P_i^{1,\ell}(t)$
are in the interval $[\ell,1]$ and we have
\begin{equation}
\label{Pi1l'}
P_i^{1,\ell}(t) = \frac{T_i^{1,0} (t,\ell)}{T_i^{1,0}(1,\ell)}=\eta_i^{1,\ell} t^i + \cdots,\quad i=0,1,\dots, k,
\end{equation}
with all leading coefficients $\eta_i^{1,\ell}>0$ and $t_{k,k}^{1,\ell}<1$. Finally, the interlacing rules
\begin{equation}\label{Interlacing}\begin{split}
t_{i,j}^{1,\ell} &\in (t_{i,j}^{1,0}, t_{i+1,j+1}^{1,0}), \ i=1,\dots, k-1, j=1,\dots, i ;\\
t_{k,j}^{1,\ell} &\in (t_{k+1,j+1}^{1,0}, t_{k,j+1}^{1,0}), \ j=1,\dots,k-1,
\end{split}
\end{equation}
hold.
\end{theorem}

\begin{remark}  As the proof below shows the condition $P_{k+1}^{1,0}(\ell)/P_{k}^{1,0}(\ell)<1$ is equivalent with $t_{k,k}^{1,\ell}<1$. In general, the orthogonal polynomial $P_k^{1,\ell} (t)$ is well defined for all $t_{k+1,1}^{1,0}<\ell<t_{k,1}^{1,0}$, but its largest root leaves the interval $[-1,1]$ and the leading coefficient becomes negative.
\end{remark}

{\it Proof.} For any polynomial $p(t)$ of degree less than $i$ we have
\begin{eqnarray*}
&& \int_{-1}^1 T_i^{1,0} (t,\ell) p(t) \, d\nu_\ell(t) \\ &=& r_i^{1,0} b_i ^{1,0}\int_{-1}^1
\left( P_{i+1}^{1,0}(t)  P_i^{1,0}(\ell)-P_{i+1}^{1,0}(\ell) P_i^{1,0}(t)   \right)p(t) \, d\chi(t)=0,
\end{eqnarray*}
and \eqref{Pi1l'} follows from the positive definiteness of the measure $d\nu_\ell(t)$ and the uniqueness of the Gram-Schmidt orthogonalization process.

We next focus on the location of the zeros of $P_i^{1,\ell}(t)$. From \eqref{kernelT} and \eqref{Pi1l'} they are solutions of the equation
\begin{equation} \label{FracEq}
\frac{P_{i+1}^{1,0}(t)}{P_i^{1,0}(t)} = \frac{P_{i+1}^{1,0}(\ell)}{P_i^{1,0}(\ell)}.
\end{equation}
For all $i<k$ the zeros of $P_{i+1}^{1,0}(t)$ and $P_i^{1,0}(t)$ are interlaced and contained in $[t_{k,1}^{1,0},t_{k,k}^{1,0}]$. Observe that ${\rm sign}\,P_i^{1,0}(\ell)=(-1)^i$, so $P_{i+1}^{1,0}(\ell)/P_i^{1,0}(\ell)<0$. The function $P_{i+1}^{1,0}(t)/P_i^{1,0}(t)$ has simple poles at $t_{i, j}^{1,0}$, $j=1,\dots, i$, and simple zeros at $t_{i+1,j}^{1,0}$, $j=1,\dots, i+1$; therefore, there is at least one solution $t_{i,j}^{1,\ell}$ of \eqref{FracEq} on every subinterval $(t_{i,j}^{1,0}, t_{i+1,j+1}^{1,0})$, $j=1,\dots,i$, which accounts for all zeros of $P_i^{1,\ell}(t)$.

When $i=k$ we note first that $P_{k+1}^{1,0}(\ell)/P_k^{1,0}(\ell)>0$. Moreover, $\ell$ is contained in the interval $(t_{k+1,1}^{1,0},t_k^{1,0})$, so we can account similarly for only the first $k-1$ solutions of \eqref{FracEq}, namely
\[t_{k,j}^{1,\ell} \in (t_{k+1,j+1}^{1,0}, t_{k,j+1}^{1,0}), \ j=1,\dots,k-1.\]
This establishes the interlacing properties \eqref{Interlacing}. To account for the last zero of $P_{k}^{1,\ell}(t)$ we utilize the fact that $P_{k+1}^{1,0}(t)/P_k^{1,0}(t)>0$ for $t\in (t_{k+1,k+1}^{1,0},\infty)$. As $\lim_{t\to \infty}P_{k+1}^{1,0}(t)/P_k^{1,0}(t) = \infty$, we have one more solution $t_{k,k}^{1,0}$ of \eqref{FracEq}.

Since $P_{k+1}^{1,0}(\ell)/P_{k}^{1,0}(\ell)<1$, we conclude that $t_{k,k}^{1,0}<1$
because $P_{k+1}^{1,0}(1)/P_k^{1,0}(1)=1$. Comparison of coefficients in \eqref{Pi1l'}  yields $\eta_k^{1,\ell}>0$.
\hfill $\Box$


\section{Construction of the Levenshtein-type polynomials}

Given some $\ell>-1$, we choose $k=k(\ell)$ to be the largest $k$ such that the condition $\ell<t_{1,k}^{1,0}$ is satisfied.

We first construct the polynomials $P_i^{1,\ell,s} (t)$ utilizing the system
$\{ P_i^{1,\ell} (t) \}_{i=0}^k$ from the previous section. The positive definiteness of the measure $d\nu_\ell(t)$ implies that
\[r_i^{1,\ell}:=\left( \int_{-1}^1 \left( P_i^{1,\ell} (t)\right)^2 \,d\nu_\ell (t) \right)^{-1} >0, \quad i=0,1,\dots,k-1.   \]
The three-term recurrence relation from Corollary \ref{cor_ortho} can be written as
\[ (t-a_i^{1,\ell}) P_i^{1,\ell} (t)=b_i^{1,\ell} P_{i+1}^{1,\ell} (t)+c_i^{1,\ell} P_{i-1}^{1,\ell} (t), \quad i=1,2, \dots ,k-1,
\]
where
\[ P_0^{1,\ell}(t)=1, \ P_1^{1,\ell}(t) = \frac{(n\ell +1)t+\ell +1}{(n+1)\ell +2}, \]
\[ b_i^{1,\ell}=\frac{\eta_{i+1}^{1,\ell}}{\eta_i^{1,\ell}}>0, \ c_i^{1,\ell}=\frac{r_{i-1}^{1,\ell} b_{i-1}^{1,\ell}}{b_i^{1,\ell}}>0, \ a_i^{1,\ell}=1-b_i^{1,\ell}-c_i^{1,\ell}. \]

By Corollary \ref{cor_ortho} we have that the zeros of $\{ P_i^{1,\ell} (t) \}$ interlace; i.e.
\[ t_{j,i}^{1,\ell}<t_{j-1,i}^{1,\ell}<t_{j,i+1}^{1,\ell}, \ \ i=1,2,\dots, j-1. \]

We next consider the Christoffel-Darboux kernel (depending on $\ell$) associated with the polynomials $\{ P_j^{1,\ell} \}_{j=0}^k$:
\begin{eqnarray}
\label{Chr_kernel_q}
R_i (x,y;\ell) &:=& \sum_{j=0}^i r_j^{1,\ell} P_j^{1,\ell}(x)  P_j^{1,\ell} (y) \\
&=& r_i^{1,\ell} b_i^{1,\ell} \frac{P_{i+1}^{1,\ell} (x)  P_i^{1,\ell} (y) - P_{i+1}^{1,\ell} (y)  P_i^{1,\ell} (x)}{x-y},
\ \ 0\leq i \leq k-1 .
\end{eqnarray}
Given $t_{k,k}^{1,0} \leq s \leq \min \{1, t_{k,k}^{1,\ell}\}$, we define
\begin{eqnarray}
\label{poly_sub}
P_{k-1}^{1,\ell,s}(t) &:=& \frac{R_{k-1} (t,s;\ell)}{R_{k-1} (1,s;\ell)}\nonumber \\
&=& \frac{1-s}{1- P_k^{1,\ell} (s)/P_{k-1}^{1,\ell} (s)}
\frac{P_k^{1,\ell} (t) -  P_{k-1}^{1,\ell} (t) P_k^{1,\ell} (s)/P_{k-1}^{1,\ell} (s) }{t-s}.
 \end{eqnarray}

We now define the Levenshtein-type polynomial
\begin{equation}\label{f_{2k}}
f_{2k}^{(n,\ell,s)}(t):= (t-\ell)(t-s)\left( P_{k-1}^{1,\ell,s} (t)\right)^2,
\end{equation}
and proceed with an investigation its properties.

\begin{theorem}\label{rootsofP1ls}
Let $n$, $\ell$, $s$ and $k$ be such that $t_{k+1,1}^{1,0}<\ell<t_{k,1}^{1,0}$, $P_{k+1}^{1,0}(\ell)/P_{k}^{1,0}(\ell)<1$,
$t_{k,k}^{1,0} \leq s \leq t_{k,k}^{1,\ell}$, and $P_k^{1,\ell} (s)/P_{k-1}^{1,\ell} (s)>P_k^{1,\ell}(\ell)/P_{k-1}^{1,\ell}(\ell)$.
Then the polynomial
$P_{k-1}^{1,\ell,s} (t)$ has $k-1$ simple zeros $\beta_1<\beta_2<\cdots<\beta_{k-1}$ such that
$\beta_1 \in (\ell,t_{k-1,1}^{1,\ell})$ and $\beta_{i+1} \in (t_{k-1,i}^{1,\ell},t_{k-1,i+1}^{1,\ell})$,
$i=1,2,\ldots,k-2$.
\end{theorem}

\medskip

{\it Proof.} The proof is similar to that of Theorem \ref{rootsofP1l}.
It follows from \eqref{poly_sub} that the roots of the
equation
\[ \frac{P_k^{1,\ell} (t)}{P_{k-1}^{1,\ell} (t)}=\frac{P_k^{1,\ell} (s)}{P_{k-1}^{1,\ell} (s)} \]
are $s$ and the zeros of $P_{k-1}^{1,\ell,s} (t)$, say $\beta_1<\beta_2<\cdots<\beta_{k-1}$.

The function $P_k^{1,\ell} (t)/P_{k-1}^{1,\ell} (t)$ has $k-1$ simple poles at the
zeros $t_{k-1,i}^{1,\ell}$, $i=1,2,\ldots,k-1$, of $P_{k-1}^{1,\ell}(t)$. Therefore,
there is a zero of $P_{k-1}^{1,\ell,s} (t)$ in each interval $(t_{k-1,i}^{1,\ell},t_{k-1,i+1}^{1,\ell})$,
$i=1,2,\ldots,k-2$, which accounts for $k-2$ zeros.

Since $P_k^{1,\ell}(s)/P_{k-1}^{1,\ell} (s)<0$ and the function $P_k^{1,\ell} (t)/P_{k-1}^{1,\ell} (t)$
increases from $-\infty$ to 1 for $t \in (t_{k-1,k-1}^{1,\ell},1]$, we have the root $s$ in this interval.
Finally, in the interval $[-\infty,t_{k-1,1}^{1,\ell})$, the function $\frac{P_k^{1,\ell} (t)}{P_{k-1}^{1,\ell} (t)}$ increases
from $-\infty$ to $+\infty$ and the condition
$P_k^{1,\ell} (s)/P_{k-1}^{1,\ell} (s)>P_k^{1,\ell}(\ell)/P_{k-1}^{1,\ell}(\ell)$ implies that
the smallest zero $\beta_1$ of $P_{k-1}^{1,\ell,s} (t)$ lies in the interval
$(\ell,t_{k-1,1}^{1,\ell})$.
\hfill $\Box$

\medskip

The next theorem is an analog of Theorem 5.39 from \cite{Lev}.
It involves the zeros of $f_{2k}^{(n,\ell,s)} (t)$ to form a right end-point Radau quadrature formula with positive weights.

 \begin{theorem} \label{QFtheorem} Let $\beta_1<\beta_2<\dots<\beta_{k-1}$
 be the zeros of the polynomial $P_{k-1}^{1,\ell,s}(t)$. Then the Radau quadrature formula
\begin{equation}
\label{QF}
\begin{split}
f_0 &= \int_{-1}^1 f(t) (1-t^2)^{\frac{n-3}{2}}\, dt \\
    &=\rho_0 f( \ell )+\sum_{i=1}^{k-1} \rho_i f(\beta_i) +\rho_k f(s)+\rho_{k+1} f(1)=:QF(f)
\end{split}
\end{equation}
is exact for all polynomials of degree at most $2k$ and has positive weights $\rho_i>0$, $i=0,1,\dots, k+1$.
\end{theorem}

\medskip

{\it Proof.} Let us denote with $L_i (t)$, $i=0,1,\dots,k+1$,
the Lagrange basic polynomials  generated by the nodes $\beta_0:=\ell<\beta_1<\dots<\beta_{k-1}<\beta_k:=s<1=:\beta_{k+1}$.
Defining $\rho_i:=\int_{-1}^1 L_i(t) d\mu(t)$, $i=0,1,\dots,k+1$,
we observe that \eqref{QF} is exact for the Lagrange basis and hence for all polynomials of degree $k+1$.

We write any polynomial $f(t)$ of degree at most $2k$ as
\[ f(t)=q(t)(t-\ell)(t-s)(1-t)P_{k-1}^{1,\ell,s} (t)+g(t), \]
where $q(t)$ is of degree at most $k-2$ and $g(t)$ is of degree at most $k+1$.
Then the orthogonality of $P_{k-1}^{1,\ell,s}(t)$ to all polynomials of degree at
most $k-2$ with respect to the measure $d\nu_{\ell,s}(t)=(t-\ell)(s-t)d\chi(t)$ and the fact
that $QF(f)=QF(g)$ show the exactness of the quadratic formula for polynomials up to degree $2k$, namely
\[ \int_{-1}^1 f(t)  d\mu(t) = \int_{-1}^1 g(t) d\mu(t) =QF(g)=QF(f).\]

We next show the positivity of the weights $\rho_i$, $i=0,\dots, k$.
Substituting in \eqref{QF} the polynomial  $f(t)=(s-t)(1-t){\left( {P_{k-1}^{1,\ell,s}}(t)\right)}^2$ of degree $2k$, we obtain
\begin{eqnarray*}
\rho_0 (s-\ell)(1-\ell)\left( P_{k-1}^{1,\ell,s} (\ell)\right)^2=f_0&=&\int_{-1}^1 (s-t)(1-t)\left( P_{k-1}^{1,\ell,s} (t)\right)^2 \, d\mu(t)\\
&=&\int_{-1}^1 \left( P_{k-1}^{1,\ell,s} (t)\right)^2 \, d\nu_s (t)>0 ,
\end{eqnarray*}
from which we derive $\rho_0 >0$.

To derive that $\rho_i>0$ for $i=1,2,\dots,k-1$, we substitute
\[ f(t)=(1-t)(t-\ell)(s-t) u_{k-1,i}^2 (t) \]
in \eqref{QF}, where $u_{k-1,i}(t)=P_{k-1}^{1,\ell,s}(t)/(t-\beta_i)$. Then clearly
\[
\rho_i (1-\beta_i)(s-\beta_i)(\beta_i-\ell)u_{k-1,i}^2 (\beta_i)=f_0=\int_{-1}^1 u_{k-1,i}^2 (t) \, d\nu_{\ell,s}(t)>0 .
\]

Similarly, utilizing the polynomial $f(t)=(1-t)(t-\ell)\left(P_{k-1}^{1,\ell,s}(t)\right)^2$ of degree $2k$ and
the positive definiteness of the measure $d\nu_\ell (t)$ up to degree $k-1$ we show that $\rho_k >0$.

Finally, we compute the weight $\rho_{k+1}$ and show that it is positive.
In this case we use $f(t)=f_{2k}^{(n,\ell,s)}(t)$ in \eqref{QF} and easily find that
\[ \rho_{k+1}=\frac{f_0}{f(1)}=\frac{f_{2k,0}^{(n,\ell,s)}}{f_{2k}^{(n,\ell,s)}(1)}=\frac{f_{2k,0}^{(n,\ell,s)}}{(1-s)(1-\ell)}. \]

Computing $f_0=f_{2k,0}^{(n,\ell,s)}$ using \eqref{poly_sub} (recall that $P_{k-1}^{1,\ell,s}(1)=1$) we get
\begin{equation}
\label{f0}
\begin{split}
f_0&=
\int_{-1}^1 (t-\ell)(s-t)(1-t)P_{k-1}^{1,\ell,s} (t) \frac{P_{k-1}^{1,\ell,s}(t)-P_{k-1}^{1,\ell,s} (1)}{t-1}  d\mu(t) \\
& \quad \quad + \int_{-1}^1 (t-\ell)(t-s)P_{k-1}^{1,\ell,s} (t)  d\mu(t) \\
&= \frac{1-s}{1-P_k^{1,\ell} (s)/P_{k-1}^{1,\ell} (s)}   \int_{-1}^1 (t-\ell)\left(  P_k^{1,\ell} (t)-\frac{P_k^{1,\ell} (s)}{P_{k-1}^{1,\ell} (s)} P_{k-1}^{1,\ell} (t) \right) \, d\mu(t).
\end{split}
\end{equation}
By Lemma \ref{Pos1} and the fact that $-P_k^{1,\ell} (s)/P_{k-1}^{1,\ell} (s) >0$ we have that the integrand in \eqref{f0} is positive definite and in particular its zero-th coefficient (which is the integral in \eqref{f0}) is positive. This proves the theorem. \hfill $\Box$

\medskip

For any fixed $-1 < \ell < t_{1,k}^{1,0}$ and $t_{k,k}^{1,0} < s < t_{k,k}^{1,\ell}$
the Levenshtein-type bound is defined to be
\[ L_{2k}(n;[\ell,s]):=\frac{1}{\rho_{k+1}}=\frac{f_{2k}^{n,\ell,s}(1)}{f_0}=\frac{(1-\ell)(1-s)}{f_0}. \]

\section{Bounding cardinalities and energies }

In the proof of the positive definiteness of his polynomials Levenshtein uses what he called the
strengthened Krein condition
\[ (t+1)P_i^{1,1}(t)P_j^{1,1}(t) \in {\mathcal F}_{>} \]
(see \cite[(3.88) and (3.92)]{Lev}). We need a following modification.

\begin{definition}
We say that
the polynomials $\{P_i^{1,\ell}(t)\}_{i=0}^k$ satisfy $\ell$-strengthened Krein condition
if
\begin{equation}
\label{LSK}
(t-\ell)P_i^{1,\ell}(t)P_j^{1,\ell}(t) \in {\mathcal F}_{>}
\end{equation}
for every $i,j \in \{0,1,\ldots,k\}$ except for $i=j=k$.
\end{definition}

The strengthened Krein condition holds true for every $i$ and $j$ by a classical result of Gasper
\cite{Gas}. However, the $\ell$-strengthened Krein condition is not true for every $\ell$, and for
fixed $\ell$, is not true for every $k$.

For fixed $n$ and $k$, denote
\[ \ell(n,k):=\sup\{ \ell \in [-1,0] : \mbox{$\ell$-strengthened Krein condition holds true} \}. \]

Our computations ensure strong evidence that the following conjecture
is true.

\begin{conjecture}
\label{conj-lsk}
For fixed $n$ and $k$ the condition \eqref{LSK} holds true for every $\ell \in [-1,\ell(n,k)]$.
\end{conjecture}

\begin{remark} The Christoffel-Darboux formula
\[ (t-\ell)P_k^{1,\ell} (t)=\frac{(1-\ell)\left( P_{k+1}^{1,0} (t)-P_k^{1,0} (t)P_{k+1}^{1,0} (\ell)/P_{k}^{1,0} (\ell) \right)}{1-P_{k+1}^{1,0} (\ell)/P_{k}^{1,0} (\ell)} \]
yields easily that the inequality $ P_{k+1}^{1,0} (\ell)/P_{k}^{1,0} (\ell)<1$
is a necessary condition for the $\ell$-strengthened condition to hold. Therefore,
we assume from now on that it holds (see the hypothesis of Theorems \ref{rootsofP1l} and \ref{rootsofP1ls}).
\end{remark}

Our computations suggest also that $\ell(n,k)$ is always less (but not much less!) than $t_{k,1}^{1,0}$ and the smallest root
of the equation $P_{k+1}^{1,0}(t)/P_{k}^{1,0}(t)=1$. Hence the $\ell$-strengthened Krein condition
is stronger than the conditions imposed so far. This corresponds to the Levenshtein's theory, where the strengthened
Krein condition appears to be the most significant obstacle.

The following Lemma demonstrates the reasonableness of Conjecture \ref{conj-lsk}.

\begin{lemma} \label{Pos1} The polynomials $\{ (t-\ell) P_i^{1,\ell} (t)\}_{i=0}^{k-1}$ are strictly positive definite
provided that $-1\leq \ell <t_{k,1}^{1,0}$.
\end{lemma}

{\it Proof.} From the definition \eqref{kernelT} of the kernels $T_i^{1,0} (x,y)$ and \eqref{Pi1l'} we have that
\begin{eqnarray*}
(t-\ell)P_i^{1,\ell} (t)&=& \frac{(1-\ell)(P_{i+1}^{1,0}(t)P_i^{1,0}(\ell)-P_i^{1,0}(t)P_{i+1}^{1,0}(\ell))}{P_{i+1}^{1,0}(1)  P_i^{1,0}(\ell)-P_i^{1,0}(1) P_{i+1}^{1,0}(\ell) }\\
&=&
\frac{1-\ell}{ 1 -P_{i+1}^{1,0}(\ell)/P_i^{1,0}(\ell)}
\left(P_{i+1}^{1,0}(t)  -\frac{P_{i+1}^{1,0}(\ell) }{P_i^{1,0}(\ell)} P_i^{1,0}(t) \right)
\end{eqnarray*}
Since the zeros of $\{ P_i^{1,0} (t) \}$ interlace we have that for all $i \leq k$
the zeros of $P_i^{1,0} (t) $ lie in the interval $[t_{k,1}^{1,0},t_{k,k}^{1,0} ]$
we have that $P_{i+1}^{1,0}(\ell)/P_i^{1,0}(\ell)<0$ for all $i\leq k-1$.
Indeed, the numerator and denominator polynomials have different signs on $(-\infty,t_{k,1}^{1,0})$.
Since $P_i^{1,0} (t) $ are strictly positive definite
(see \cite[Eq. (3.91)]{Lev}), we conclude the proof of the Lemma. \hfill
$\Box$

\medskip


\begin{lemma}
\label{Pos2} The polynomials $(t-\ell)(t-s)P_{k-1}^{1,\ell,s} (t)$ and
$(t-\ell)P_{k-1}^{1,\ell,s} (t)$ are strictly positive definite.
\end{lemma}

{\it Proof.}
The interlacing property of $\{P_i^{1,\ell}\}$ implies that $ \frac{P_k^{1,\ell} (s)}{P_{k-1}^{1,\ell} (s)}<0$.
Applying Lemma \ref{Pos1} one concludes that $(t-\ell)(t-s)P_{k-1}^{1,\ell,s}(t)$
is positive definite. Furthermore, as  $P_j^{1,\ell}(s)>0$ and $(t-\ell)P_j^{1,\ell} (t)$
is strictly positive definite for $j=0,1,\dots,k-1$, we derive using \eqref{Chr_kernel_q} that $(t-\ell)P_{k-1}^{1,\ell,s} (t)$ is also strictly positive definite.
\hfill $\Box$

\medskip

In our main result we use the $\ell$-strengthened Krein condition relying on the following observation. For fixed $n$, $k$,
and $\ell$, we check numerically whether \eqref{LSK} is satisfied for every
pair $(i,j)$, $i,j \in \{0,1,\ldots,k\}$, except for $i=j=k$. This is done for every $\ell=-1+m\varepsilon$,
$\varepsilon=10^{-3}$ (of course, doing the step $\varepsilon$ smaller is only a matter of computations), $m=1,2,\ldots$,
until \eqref{LSK} holds true. In practice, when one needs to compute bounds in the class $C(\ell,s)$, he can consider
instead $C(\ell_0,s)$, where $\ell_0<\ell$ is the largest
for which the $\ell$-strengthened Krein condition holds true.

The next assertion is the analog of Theorem 5.42 of \cite{Lev}. It uses a seemingly
weaker\footnote{In fact, we suspect that both
conditions are equivalent.} version of the
$\ell$-strengthened Krein condition.

 \begin{theorem} \label{PD_Lev}
Let $n$, $k$, and $\ell$ be such that the polynomials $(t-\ell)P_i^{1,\ell}(t)P_j^{1,\ell}(t)$ are positive definite for
$i \in \{k,k-1\}$ and every $j \leq k-1$. Let $t_{k,k}^{1,0} \leq s \leq t_{k,k}^{1,\ell}$ and $P_k^{1,\ell} (s)/P_{k-1}^{1,\ell} (s)>P_k^{1,\ell}(\ell)/P_{k-1}^{1,\ell}(\ell)$. Then the Levenshtein-type polynomial $f_{2k}^{(n,\ell,s)}(t)$ is positive definite.
 \end{theorem}

\medskip

{\it Proof.} It follows from the definition \eqref{poly_sub} that the Levenshtein-type polynomial can be represented as
follows
\begin{equation}
\label{f_lev_pd1}
f_{2k}^{(n,\ell,s)}(t)= c(t-\ell)\left(P_k^{1,\ell} (t)+c_1 P_{k-1}^{1,\ell} (t)\right)
\sum_{i=0}^{k-1} r_i^{1,\ell} P_i^{1,\ell} (t)P_i^{1,\ell}(s),
\end{equation}
where $c=(1-s)/(1+c_1)R_{k-1}(1,\ell,s)>0$ and $c_1=-P_k^{1,\ell} (s)/P_{k-1}^{1,\ell} (s)>0$ under
the assumptions for $\ell$ and $s$. Since $P_i^{1,\ell} (s)>0$ for $0 \leq i \leq k-1$, the polynomial $f_{2k}^{(n,\ell,s)}(t)$
becomes positive linear combination of terms like
$(t-\ell)P_i^{1,\ell}(t)P_j^{1,\ell}(t)$, where $i \in \{k,k-1\}$ and $j \leq k-1$.  \hfill $\Box$

\medskip

The main result in this paper is the following.

\begin{theorem}
\label{main-th-subintervals}
Assume that $\ell \in [-1,t_{k,1}^{1,0})$ and $s \in (t_{k,k}^{1,0},t_{k,k}^{1,\ell})$ and that
the $\ell$-strengthened Krein condition holds true.
Then
\begin{equation}
\label{L-like-bound}
 \mathcal{A}(n;[\ell,s]) \leq \frac{f_{2k}^{(n,\ell,s)}(1)}{f_0}=\frac{1}{\rho_{k+1}}.
\end{equation}
Furthermore, for $h$ being an absolutely monotone function, and for $M$ determined
by $f_{2k}^{(n,\ell,s)}(1)=Mf_0$,
the Hermite interpolant\footnote{The notation $g=H(f;h)$ is taken from \cite{CK}; it signifies that
$g$ is the Hermite interpolant to the function $h$ at the zeros (taken with their multiplicity) of $f$.}
\[ g(t)=H((t-s)f_{2k}^{(n,\ell,s)}(t);h) \]
belongs to $\mathcal{G}_{n,\ell;h}$, and, therefore,
\begin{equation}
\label{ULB-like-bound}
 \mathcal{E}(n,M,\ell;h) \geq M(Mg_0-g(1))=M^2\sum_{i=0}^{k} \rho_i h(\beta_i).
\end{equation}
\end{theorem}

{\it Proof.} We first verify the positive definiteness of the polynomials
$f_{2k}^{(n,\ell,s)}(t)$  and $g(t)$. We have $f_{2k}^{(n,\ell,s)}(t) \in \mathcal{F}_>$ by Theorem \ref{PD_Lev}.

Denote by $t_1\leq t_2\leq \cdots \leq t_{2k}$ the zeros of $f_{2k}^{(n,\ell,s)}(t)$ counting multiplicity. Observe, that $t_1=\ell$, $t_{2i}=t_{2i+1}=\beta_i$, $i=1,\dots,k-1$, and $t_{2k}=s$.
It follows from \cite[Lemma 10]{CW} that the polynomial
\[ g(t):=H((t-s)f_{2k}^{(n,\ell,s)}(t);h) \]
is a linear combination with nonnegative coefficients
of the partial products
\[ \prod_{j=1}^{m} (t-t_j), \ \ m=1, 2,\ldots,2k. \]

Since $t_{2i}$, $i=1,\ldots,k$, are the roots of $P_k^{1,\ell}(t)+\alpha P_{k-1}^{1,\ell}(t)$
(see \eqref{poly_sub}) it follows from \cite[Theorem 3.1]{CK} that the partial products
$\prod_{j=1}^{m} (t-t_{2j})$, $m=1,\ldots,k-1$, have positive coefficients when expanded in terms
of the polynomials $P_{i}^{1,\ell}(t)$. Then $g(t)$ is a linear combination with positive coefficients of
terms $(t-\ell)P_i^{1,\ell}(t)P_j^{1,\ell}(t)$ and the
last partial product which is in fact $f_{2k}^{(n,\ell,s)}(t)$.
Now the positive definiteness of $g(t)$ follows from the validity of the $\ell$-strengthened Krein condition
and Theorem \ref{PD_Lev}.

The bounds \eqref{L-like-bound} and \eqref{ULB-like-bound} now hold true since $f_{2k}^{(n,\ell,s)}(t) \leq 0 $
for every $t \in [\ell,s]$ (see \eqref{f_{2k}}) and $g(t) \leq h(t)$ for every $t \in [\ell,1)$ by \cite[Lemma 9]{CW}.

The expressions of the bounds via the weights $\rho_i$ and the nodes $\beta_i$ follow from Theorem
\ref{QFtheorem}.
\hfill $\Box$

\section{Examples and numerical results}

\subsection{On the $\ell$-strengthened Krein condition}

In the table below we present our computations of the value of $\ell(n,k)$, the maximum $\ell$ for fixed $n$ and
$k$, such that the $\ell$-strengthened Krein condition is true.

\begin{center}
{\bf Table.} Conjectured values of $\ell(n,k)$ for $3 \leq n,k \leq 10$. The rows after the
corresponding $\ell(n,k)$ show the value of the smallest root of the equation $P_{k+1}^{1,0}(t)/P_{k}^{1,0}(t)=1$.
The real numbers are truncated after the third digit.

 \label{tab:1}

\medskip

\begin{tabular}{|c|c|c|c|c|c|c|c|c|c|c|}
  \hline
$k$ & 10 & 9 & 8 & 7 & 6 & 5 & 4 & 3 & 2 & 1   \\ \hline
$\ell(3,k)$ & $-.979$ & $-.974$ & $-.969$ & $-.962$ & $-.951$ & $-.936$ & $-.912$ & $-.870$ & $-.787$ & $-.577$   \\ \hline
& $-.978$ & $-.973$ & $-.967$ & $-.959$ & $-.948$ & $-.930$ & $-.902$ & $-.854$ & $-.754$ & $-.500$   \\ \hline
$\ell(4,k)$ & $-.967$ & $-.961$ & $-.953$ & $-.942$ & $-.927$ & $-.906$ & $-.874$ & $-.821$ & $-.723$ & $-.499$   \\ \hline
& $-.965$ & $-.958$ & $-.950$ & $-.938$ & $-.922$ & $-.897$ & $-.860$ & $-.796$ & $-.676$ & $-.400$   \\ \hline
$\ell(5,k)$ & $-.955$ & $-.947$ & $-.936$ & $-.923$ & $-.905$ & $-.879$ & $-.840$ & $-.779$ & $-.672$ & $-.447$   \\ \hline
& $-.952$ & $-.944$ & $-.932$ & $-.917$ & $-.896$ & $-.866$ & $-.821$ & $-.748$ & $-.615$ & $-.333$   \\ \hline
$\ell(6,k)$ & $-.942$ & $-.933$ & $-.921$ & $-.904$ & $-.883$ & $-.853$ & $-.810$ & $-.744$ & $-.631$ & $-.408$   \\ \hline
& $-.939$ & $-.929$ & $-.915$ & $-.897$ & $-.872$ & $-.838$ & $-.787$ & $-.706$ & $-.566$ & $-.285$   \\ \hline
$\ell(7,k)$ & $-.930$ & $-.919$ & $-.905$ & $-.887$ & $-.863$ & $-.830$ & $-.783$ & $-.712$ & $-.597$ & $-.377$   \\ \hline
& $-.926$ & $-.914$ & $-.898$ & $-.878$ & $-.850$ & $-.811$ & $-.755$ & $-.670$ & $-.526$ & $-.250$   \\ \hline
$\ell(8,k)$ & $-.918$ & $-.906$ & $-.890$ & $-.870$ & $-.843$ & $-.808$ & $-.758$ & $-.685$ & $-.568$ & $-.353$   \\ \hline
& $-.914$ & $-.900$ & $-.882$ & $-.859$ & $-.828$ & $-.787$ & $-.727$ & $-.638$ & $-.492$ & $-.222$   \\ \hline
$\ell(9,k)$ & $-.907$ & $-.893$ & $-.876$ & $-.854$ & $-.825$ & $-.788$ & $-.736$ & $-.660$ & $-.543$ & $-.333$   \\ \hline
& $-.901$ & $-.886$ & $-.866$ & $-.841$ & $-.808$ & $-.764$ & $-.702$ & $-.610$ & $-.463$ & $-.200$   \\ \hline
$\ell(10,k)$ & $-.895$ & $-.880$ & $-.862$ & $-.838$ & $-.808$ & $-.769$ & $-.715$ & $-.638$ & $-.520$ & $-.316$   \\ \hline
& $-.889$ & $-.872$ & $-.851$ & $-.824$ & $-.789$ & $-.743$ & $-.678$ & $-.585$ & $-.439$ & $-.181$   \\ \hline
  \end{tabular}
	
\end{center}

\subsection{System of bounds for fixed $n$ and $\ell$}

We present here as example the system of bounds for $\mathcal{A}(n;[\ell,s])$, where $n=4$ and $\ell=-0.95$ are
fixed and $s$ is varying. According to the above table, the $\ell$-strengthened Krein conditions holds true for $k \leq 7$ and
corresponding bounds
\[ \mathcal{A}(4;[-0.95,s]) \leq L_{2k}(4;[-0.95,s])=1/\rho_{k+1}, \ k=1,2,\ldots,7, \]
hold true.

On the figure below we show the first four bounds
\[ L_{2k}(4;[-0.95,s]), \ k=1,2,3,4, \]
together with the Levenshtein odd degree bounds $L_{2u-1}(4,s)$, $u=1,2,3,4$. The subscripts are missed for short.
The behaviour of the bounds is as follows. For $s \in [-0.95,t_1^{1,0}]$, $t_1^{1,0}=-1/4$, the Levenshtein bound $L_1(4,s)$ is
better, then for $s \in [t_1^{1,0},0.0175]$ our bound $L_{2}(4;[-0.95,s])=1/\rho_{2}$ is better,
for $s \in [-0.0175,t_2^{1,0}]$, $t_2^{1,0} \approx 0.27429$, the Levenshtein bound $L_3(4,s)$ is
better, then for $s \in [t_2^{1,0},0.4195]$ our bound $L_{4}(4;[-0.95,s])=1/\rho_{4}$ is better, etc.
This is the typical situation for all reasonable values of $n$ and $\ell$ we have checked.

\begin{center}	
   \begin{figure}[h!]	
	\includegraphics[width=1.0\textwidth]{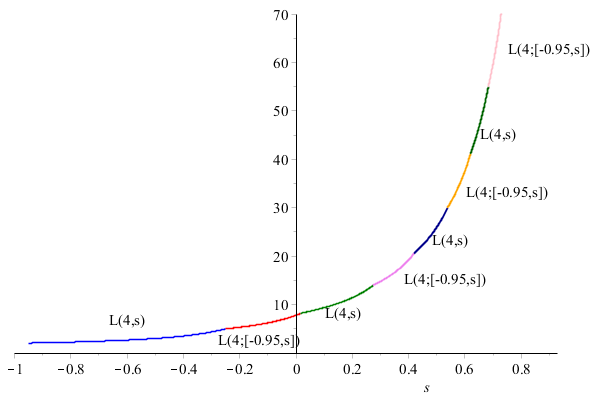}
	\end{figure}
\end{center}

\subsection{Bounds for $\mathcal{E}(n,M,\ell;h)$}

We use the system of bounds from Section 5.2 to derive our ULB-like bounds for $\mathcal{E}(n,M,\ell;h)$.
Given $n$, $\ell$, and $M$ we consecutively construct the polynomials $f_{2k}^{n,\ell,s}(t)$ and their
bounds as above until we reach the maximum $k$ such that $f_{2k}^{n,\ell,s}(t) \in \mathcal{F}_>$ and
the equality
\[ f_{2k}^{n,\ell,s}(1)=Mf_0 \]
holds. Then, as Theorem \ref{main-th-subintervals} states, we construct the interpolant
\[ g(t)=H((t-s)f_{2k}^{(n,\ell,s)}(t);h) \]
and compute the bound
\[ \mathcal{E}(n,M,\ell;h) \geq M(Mg_0-g(1))=M^2\sum_{i=0}^{k} \rho_i h(\beta_i). \]

\subsection{System of bounds for fixed $\ell$ and $s$}

For the case $k=2$, the Levenshtein polynomial  \eqref{f_{2k}} is given by
\[ f_{4}^{(n,\ell,s)}(t):= (t-\ell)(t-s)\left( P_{1}^{1,\ell,s} (t)\right)^2, \]
where the zero of $P_{1}^{1,\ell,s} (t)$ is $\alpha=-\frac{3+(n+2)(\ell s+\ell+s)}{(n+2)(n\ell s+\ell+s+1)}$.
Thus, from Theorem \ref{main-th-subintervals} we obtain
\begin{equation}
\label{U4}
\mathcal{A}(n;[\ell,s]) \leq \frac{f_{4}^{(n,\ell,s)}(1)}{f_0}=
\frac{n(1-\ell)(1-s)[3+(n+2)(n\ell s+\ell s+2\ell+2s+1)]}{(n+2)[n\ell^2s^2-(\ell-s)^2]-6\ell s+3}
\end{equation}
subject to
 \[ \ell+s+2\alpha \leq 0, \]
\[ \alpha^2+2(\ell+s)\alpha+\ell s+\frac{6}{n+4} \geq 0, \]
\[ (\ell+s)\alpha^2+2\alpha(\ell s+\frac{3}{n+2})+\frac{3(\ell+s)}{n+2} \leq 0.\]
The bound \eqref{U4} is attained by codes of parameters
\[ (n,M,s)=\left(3m^2-5,\frac{m^4(3m^2-5)}{2},\frac{1}{m+1}\right), \]
known only for $m=2$ (here $\ell=-1$) and 3 (here $\ell=-1/4$). Such codes are
derived from corresponding tight spherical 7-designs in dimensions $3m^2-4$ (see \cite{DGS}).

Let $n$, $M$, and $\ell$ be such that $k=2$ be the maximal value of $k$ such that $f_{2k}^{n,\ell,s}(1)=Mf_0$
holds true and the above $f_{4}^{(n,\ell,s)}(t)$ is positive definite (this fixes $s$ as well).
Then, according to Theorem \ref{main-th-subintervals}, the $h$-energy (for any absolutely monotone $h$)
bound \eqref{ULB-like-bound} is given by the polynomial $g_4(t) \in \mathcal{G}_{n,\ell;h}$
of degree 4 which interpolates $h$ by
\[ g_4(\ell)=h(\ell), \ g_4(\alpha)=h(\alpha), \ g_4^\prime(\alpha)=h^\prime(\alpha), \
g_4(s)=h(s), \ g_4^\prime(s)=h^\prime(s). \]

In right ranges for $\ell$ and $s$ both bounds are optimal in the sense that they can not be improved by using linear programming
with polynomials of degree at most 4.

{\bf Acknowledgement.} The authors thank Konstantin Delchev, Tom Hanson, and Nikola Sekulov for their computational work that independently verified Conjecture \ref{conj-lsk}.


\begin{thebibliography}{99}

\bibitem{BCV} A. Bultheel, R. Cruz-Barroso and M. Van Barel, On Gauss-type quadrature
formulas with prescribed nodes anywhere on the real line, {\it Calcolo} 47 (2010), 21-48.

\bibitem{BBMQ} B.\,Beckermann, J.\,Bustamante, R. Martinez-Cruz, and J. Quesada,
Gaussian, Lobatto and Radau positive quadrature rules with a prescribed abscissa, {\em Calcolo} 51, (2014), 319-328.

\bibitem{BHS}
S.\,Borodachov, D.\,Hardin, and E.\,Saff,
\emph{Minimal Discrete Energy on Rectifiable Sets}, Springer, 2018 (to appear).

\bibitem{BDHSS-CA} P.\,Boyvalenkov, P.\,Dragnev, D.\,Hardin, E.\,Saff, and M.\,Stoyanova.
Universal lower bounds for potential energy of spherical codes, {\it Constr. Approx.} 44, 2016, 385-415.

\bibitem{BDHSS-DCC} P.\,Boyvalenkov, P.\,Dragnev, D.\,Hardin, E.\,Saff, and M.\,Stoyanova,
Energy bounds for codes and designs in Hamming spaces, {\it Designs, Codes and Cryptography},
82(1), (2017), 411-433 (arxiv:1510.03406).

\bibitem{CK} H.\,Cohn and A.\,Kumar,
Universally optimal distribution of points on spheres,
\emph{J. Amer. Math. Soc.} 20, (2007), 99-148.

\bibitem{CW} H.\,Cohn and J.\,Woo, Three point bounds for energy minimization,
\emph{J. Amer. Math. Soc.} 25, (2012), 929-958.


\bibitem{CZ} H.\,Cohn and Y.\,Zhao,
Energy-minimizing error-correcting codes, \emph{IEEE Trans. Inform. Theory}
60, (2014), 7442-7450 (arXiv:1212.1913).


\bibitem{DR} P. J. Davis and P. Rabinowitz, \emph{Methods of Numerical Integration}, 2nd ed. Academic Press, New York (1984).

\bibitem{DGS} P.\,Delsarte, J.-M.\,Goethals, and J.\, J.\,Seidel,
Spherical codes and designs,
\emph{Geom. Dedicata} 6, (1977), 363-388.

\bibitem{Gas} G. Gasper, Linearization of the product of Jacobi polynomials, II, \emph{Canad. J.
Math.} 22, (1970), 582-593.

\bibitem{KL} G. A. Kabatyanskii and V. I. Levenshtein, Bounds for packings
on a sphere and in space, {\it Probl. Inform. Transm.} 14,
(1989), 1-17.

\bibitem{Lev92} V.\,I.\,Levenshtein, Designs as maximum codes in polynomial
metric spaces, \emph{Acta Appl. Math.} 25, (1992), 1-82.

\bibitem{Lev}
V.\,I.\,Levenshtein,
Universal bounds for codes and designs, \emph{Handbook of Coding Theory},
V.\,S.~Pless and W.\,C.~Huffman, Eds., Elsevier, 1998, Ch.~6, 499--648.

\bibitem{Sid} {  V.\,M.\,Sidel'nikov}, On extremal polynomials used to
estimate the size of codes, {\it Probl. Inform. Transm.} 16 (1980), 174-186.

\bibitem{Sze}
G.\,Szeg\H{o}, \emph{Orthogonal polynomials},
Amer. Math. Soc. Col. Publ., {\bf 23}, Providence, RI, 1939.

\bibitem{Y} V.\,A.\,Yudin, Minimal potential energy of a point system of charges,
\emph{Discret. Mat.} 4, (1992), 115-121 (in Russian); English translation: \emph{Discr. Math. Appl.} 3, (1993), 75-81.

\end{thebibliography}
\end{document}